\documentclass[10pt]{article}
\usepackage{a4wide}
\usepackage{amsmath,amssymb}
\newtheorem{theorem}{Theorem}[section]
\newtheorem{lemma}{Lemma}[section]
\newtheorem{remark}{Remark}[section]

\date{}
\begin{document}
\title {Travelling wavefronts in nonlocal diffusion equations with nonlocal delay effects}
\author {Shangjiang Guo\thanks{E-mail: shangjguo@hnu.edu.cn}\\
  College of Mathematics and Econometrics, Hunan
  University\\
  Changsha, Hunan 410082, People's Republic of China\\
  Johannes Zimmer
  \\
  Department of Mathematical Sciences, University of Bath, \\
  Bath BA2 7AY, United Kingdom}
 \maketitle
\begin{abstract}
  This paper deals with the existence, monotonicity, uniqueness and   asymptotic behaviour of travelling wavefronts for a class of temporally delayed, spatially nonlocal diffusion equations.

  {\bf Keywords.}  Travelling waves, time delay, Fisher-KPP equation

  {\bf AMS subject classifications.} 35K57, 34K05, 92D25.
\end{abstract}

\section{Introduction}
\label{sec:Introduction}

Travelling wavefront solutions play an important role in the description of the long-term behaviour of solutions to initial value problems in reaction-diffusion equations, both in the spatially continuous case and in spatially discrete situations. Such solutions are also of interest in their own right, for example to understand transitions between different states of a physical system, propagation of patterns, and domain invasion of species in population biology (see, e.g.,~\cite{Bates1,Bates2,Cahn,Chow,So}). In this paper, we study the existence, uniqueness and asymptotic stability of travelling wavefronts of the equation:
\begin{equation}
  \label{eq}
  u_t(x,t)=p u_{xx}(x,t)+d(J*u-u)(x,t)+f
  \left(u(x,t),(h**u)(x,t)\right),
\end{equation}
where $x\in\mathbb{R}$, $d\geq 0$, $p\geq 0$, and
\begin{equation*}
\aligned
  (J*u)(x,t):=&\int_{\mathbb{R}}J(x-y)u(y,t)dy,\\
  (h**u)(x,t):=&\int^{t}_{-\infty}\int_{\mathbb{R}}h(x-y,t-s)u(y,s)dyds.
\endaligned
\end{equation*}
Equation (\ref{eq}) mixes a continuous Laplacian with a nonlocal diffusion $d(J*u-u)(x,t)$, which describes that the diffusion of
density $u$ at a point $x$ and time $t$ depends not only on $u(x, t)$ but also on all the values of $u$ in a neighbourhood of $x$
through the convolution term $J*u$. In population dynamics, the reaction term $f(u, h**u)$ is usually used to describe the
recruits of population, and $h**u$ represents a weighted average of the population density both in past time and
space~\cite{Britton1,Britton2}. The nonlinear functions $f(u,v)$ and $h(u)$ satisfy the following hypotheses:
\begin{description}
\item[(F1)] $f\in C([0,K]\times [0,K],\mathbb{R})$, $f(0,0)=f(K,K)=0$,   $f(u,h**u)>0$ for all $u\in (0,K)$, $\partial_2f(u,v)\geq 0$ for all   $(u,v)\in [0,K]\times [0,K]$, where $K$ is a positive constant.
\item[(F2)] There exist some $M>0$ and $\sigma\in (0,1]$ such that
   $0\leq \partial_1f(0,0)u+\partial_2f(0,0)v-f(u,v)\leq         M(u+v)^{1+\sigma}$
   for all $(u,v)\in [0,K]\times [0,K]$ and    $\partial_1f(K,K)+\partial_2f(K,K)<0$.
 \item[(H1)] Both $J(\cdot)$ and $h (\cdot,t)$ are nonnegative, even, integrable and satisfies $\int_{\mathbb{R}}J(x)dx=1$ and $\int_0^{\infty}\int_{\mathbb{R}}h(x,t)dxdt=1$.
 \item[(H2)] There exists some $\lambda_0>0$ (possibly equal to    $\infty$) such that $\int^{\infty}_0J(x)\exp\{\lambda x\}dx<\infty$ and $\int_0^{\infty}\int_{0}^{\infty}h(x,t)\exp\{\lambda (x-ct)\}dxdt<\infty$ for all $c\geq 0$ and $\lambda\in [0,\lambda_0)$.
\end{description}

Assumptions (F1) and (F2) are standard. From (F1) we can see that (\ref{eq}) has two equilibria 0 and $K$. Furthermore, condition
(F2) together with (F1) and (H1) implies that
$\partial_1f(0,0)+\partial_2f(0,0)\geq \frac{2}{K}f(\frac{K}{2},\frac{K}{2})=\frac{2}{K}f(\frac{K}{2},h**(\frac{K}{2}))>0$, hence
0 is unstable and $K$ is stable. In this article, we will not require that $\partial_2f(0,0)> 0$.

Since Equation (\ref{eq}) involves a general diffusion kernel and delayed nonlinearity, it can be reduced to some well-known
equations if $J$, $h$, and $f$ are chosen to take a a some special form (see, for
example,~\cite{AL-OmariGourley2005,Faria2006,Schaaf1987,Wang-Li-Ruan,WengWu2008,WuZou2001,ZhaoXiao2006}).  In particular, special
cases of~\eqref{eq} include a host-vector disease model, a nonlocal population model with age structure, and a nonlocal
Nicholson's blowflies model with delay; these cases are discussed in a second paper investigating the stability of the
system~\cite{Guo:11b}. For example, choosing $d=0$, $f(u,v)=-\tau u+\tau\beta ve^{-v}$, equation (\ref{eq}) can be reduced to the
following Nicholson's Blowflies equation with spatio-temporal delays
\begin{equation*}
   u_t(x,t)=p u_{xx}(x,t)-\tau u(x,t)+\tau\beta(h**u)(x,t)\exp\{-(h**u)(x,t)\},
\end{equation*}
which was studied by Li, Ruan, and Wang~\cite{Li-Ruan-Wang2007}, and Lin~\cite{Lin2009}.
If $p=0$ and $J(x)=\frac{1}{2}[\delta(-1)+\delta(1)]$ and $h(x,t)=k(x)\delta(t-\tau)$, where $\delta(\cdot)$ is the Dirac delta function, then (\ref{eq}) reduces to the discrete reaction-diffusion equation
\begin{equation}
  \label{eq000}
  u_t(x,t)= d\cdot\Delta_1u(x,t)+f
  \left(u(x,t),(k*u)(x,t-\tau)\right),
\end{equation}
where $\Delta_1u(x,t)=\frac{1}{2}[u(x+1,t)-2u(x,t)+u(x-1,t)]$. If $f(u,v)=-au+b(v)$ and $k(x)=\delta(x)$, then (\ref{eq000}) reduces to the local equation
\begin{equation}
  \label{ex0}
  u_t(x,t)=d\cdot\Delta_1u(x,t)-au(x,t)+b(u(x,t-\tau)),\quad x\in
  \mathbb{R},\,t\geq 0,
\end{equation}
where $b\in C^1([0,\infty],\mathbb{R})$.  If $f(u,v)=g(u)$ and $g(u)$ denotes a Lipschitz continuous function satisfying
$g(u)>0=g(0)=g(1)$ for all $u\in (0,1)$, equation (\ref{eq000}) becomes
\begin{equation}
  \label{eqc2}
  u_t(x,t)= d\cdot\Delta_1u(x,t)+g(u(x,t)).
\end{equation}
On the other hand, when $d=0$ and $h(x,t) =k(x)\delta(t-\tau)$, (\ref{eq}) reduces to the following reaction–diffusion equation with discrete time delay
\begin{equation}\label{eqcon1}
  u_t(x,t)= pu_{xx}(x,t)+f
  \left(u(x,t),(k*u)(x,t-\tau)\right).
\end{equation}
Moreover, the equation (\ref{eq000}) is a spatially discrete version of (\ref{eqcon1}) with $p$ replaced by $d$. In recently
years, spatially non-local differential equations such as (\ref{eqcon1}) have attracted significant attention (see,
e.g.,~\cite{Fang,Gourley,Ou,So,Thieme, Wang-Li-Ruan}). Under some monostable assumption, Wang \textit{et al}.~\cite{Wang-Li-Ruan}
investigated the existence, uniqueness, and global asymptotical stability of travelling wave fronts. We also refer to So \emph{et
  al}.~\cite{So} for more details and some specific forms of $f$, obtained from integrating along characteristics of a structured
population model, an idea from the work of Smith and Thieme~\cite{Smith}. See also~\cite{So} for a similar model
and~\cite{Gourley0} for a survey on the history and the current status of the study of reaction diffusion equations with non-local
delayed interactions.  In particular, when $f(u,v)= v(1-u)$ and $k(u)=\delta(u)$, the equation (\ref{eqcon1}) is delayed
\emph{Fisher's equation}~\cite{Fisher} or \emph{KPP equation}~\cite{Kolmogorov}, which arises in the study of gene development or
population dynamics. When $f(u,v)=-au+b(v)$ and $k(u)=\delta(u)$, the equation (\ref{eqcon1}) is the local \emph{Nicholson's
  blowflies equation} and has been investigated in~\cite{Gourley,Gurney,Lin, Mei1}. When $f(u,v)=-au+b(1-u)v$, equation
(\ref{eqcon1}) is called the \emph{vector disease model} as proposed by Ruan and Xiao~\cite{Ruan}. When
$f(u,v)=bv\exp\{-\gamma\tau\}-\delta u^2$ and $k(u)=\frac{1}{\sqrt{4\pi \alpha\tau}}\exp\{\frac{-y^2}{4\alpha\tau}\}$, equation
(\ref{eqcon1}) is the age-structured reaction diffusion model of a single species proposed by Al-Omari \&
Gourley~\cite{Al}. Existence and stability of travelling wavefronts for the reaction-diffusion equation (\ref{eqcon1}) and its
special forms has been extensively studied in the literature.

We are interested in wave propagation phenomena. In particular, we are interested in monotone travelling waves $u(x, t) = \phi(x + ct)$ for (\ref{eq}), with $\phi$ saturating at 0 and $K$. We call $c$ the \emph{travelling wave speed} and $\phi$ the \emph{profile} of the wavefront. In order to address these questions, we need to find an increasing function $\phi(\xi)$, where $\xi= x + ct$, which is a solution of the associated travelling wave equation
\begin{equation}
  \label{eq3}
  %\left\{
  \begin{split}
      -c\phi'(\xi)+p\phi''(\xi)+ d(J*\phi-\phi)(\xi)+f(\phi(\xi),
      (h**\phi)(\xi))=0,\quad \xi\in \mathbb{R},\\
      \lim\limits_{\xi\to-\infty}\phi(\xi)=0,\qquad
      \lim\limits_{\xi\to \infty}\phi(\xi)=K,\\
      0\leq \phi(\xi)\leq K,\quad \xi\in \mathbb{R},
    \end{split} %\right.
\end{equation}
where
$$
\aligned
(J*\phi)(\xi)=&\int_{\mathbb{R}}J(y)\phi(\xi-y)dy,\\
(h**\phi)(\xi)=&\int^{\infty}_{0}\int_{\mathbb{R}}
h(y,s)\phi(\xi-y-cs)dyds.
 \endaligned
$$
 For convenience, we write $\phi(-\infty)$ and $\phi(\infty)$ as abbreviations for $\lim_{\xi\to-\infty}\phi(\xi)$ and $\lim_{\xi\to   \infty}\phi(\xi)$, respectively. Travelling wavefronts of (\ref{ex0}) have been intensively studied in recent years, see, e.g., \cite{Cahn,Carr,Chen0,Chen1,Chen2,Chow,Hsu,Hudson,Keener,MaLiao,Ma,Mallet,WuZou1,Zinner1991,Zinner1992,Zinner}. Zinner \emph{et al.}~\cite{Zinner} addressed the existence and minimal speed of travelling wavefront for (\ref{eqc2}).  Recently, based on~\cite{Chen1,Chen2}, Chen \emph{et al}.~\cite{Chen0} investigated the uniqueness and asymptotic behaviour of travelling waves for (\ref{eq000}) with $d=2$. To the best of our knowledge, however, there are no results regarding the existence, uniqueness, monotonicity, asymptotic behaviour, and asymptotic stability of travelling waves for an equation as general as (\ref{eq}).

There is an enormous amount of work on related equations which is impossible even to sketch. We only mention the work of Coville and coworkers, where the nonlinearity is local, but general nonlocal expressions instead of the nonlocal diffusion equation are considered (e.g.,~\cite{Coville}). Some methods are similar, such as the use of super- and subsolutions. For interesting work on a Fisher-KPP equation with a non-local saturation effect, where no maximum principle holds, we refer to~\cite{Berestycky}.

We shall establish the existence, uniqueness, monotonicity, asymptotic behaviour of travelling waves for (\ref{eq}) under the assumptions (F1), (F2), (H1), and (H2).

\begin{theorem}\label{thm-existence} Under assumptions (F1), (F2),     (H1), and (H2), there exists a minimal wave speed $c^*>0$ such that   for each $c\geq c^*$, equation (\ref{eq}) has a travelling wavefront   $\phi(x+ct)$ satisfying (\ref{eq3}).
       Moreover,
  \begin{enumerate}
  \item the solution $\phi$ of (\ref{eq3}) is unique up to a     translation.
  \item Every solution $\phi$ of (\ref{eq3}) is strictly monotone,     i.e., $\phi'(\xi)>0$ for all $\xi\in \mathbb{R}$.
  \item Every solution $\phi$ of (\ref{eq3}) satisfies $0 <\phi(\cdot)     <K$ on $\mathbb{R}$.
  \item Any solution of (\ref{eq3}) satisfies     $\lim_{\xi\to-\infty}\phi'(\xi)/\phi(\xi)=\lambda$, with $\lambda$     being the minimal positive root of
    \begin{equation}\label{1eq}
      c\lambda-p\lambda^2-d\left[H(\lambda)-1\right]-\partial_1f(0,0)
      -\partial_2f(0,0)G(c,\lambda)=0,
    \end{equation}
    where
    \begin{equation*}
      \aligned
      H(\lambda)=&\int_{\mathbb{R}}J(y)\exp\{-\lambda y\}dy,\\
      G(c,\lambda)=&\int^{\infty}_0\int_{\mathbb{R}}h(y,s)\exp\{-\lambda (y+cs)\}dyds.
      \endaligned
    \end{equation*}
    for $\lambda\in \mathbb{C}$ with $\mathrm{Re}\lambda<\lambda_0$.
  \item Any solution of (\ref{eq3}) satisfies $\lim_{\xi\to \infty}\phi'(\xi)/[K-\phi(\xi)]=\gamma$, with $\gamma$ being the
    unique positive root of
    \begin{equation}\label{2eq}
      c\gamma+p\gamma^2+d\left[H(-\gamma)-1\right]+\partial_1f(K,K)
      +\partial_2f(K,K)G(c,-\gamma)=0.
    \end{equation}
  \end{enumerate}
\end{theorem}

We remark that if $c>c^*$, equation (\ref{1eq}) has exactly two real roots, both positive.

This paper is organised as follows. In Section~\ref{sec:Notat-auxil-results}, we provide some preliminary results; in Section~\ref{sec:Existence}, we establish the existence of a travelling wavefront, using the monotone iteration method developed by Wu and Zou~\cite{WuZou1} with a pair of super- and sub-solutions. In particular, Theorems~\ref{thm-ex} and \ref{thm-cri} establish the existence part of Theorem~\ref{thm-existence}. To derive the monotonicity and uniqueness of wave profiles (Section~\ref{sec:Monot-Uniq}), we shall first apply Ikehara's theorem in Section~\ref{sec:Asymptotic-behaviour} to study the asymptotic behaviour of wave profiles. This idea originated in Carr and Chmaj's paper~\cite{Carr}, where the authors study the uniqueness of waves for a nonlocal monostable equation. Theorem~\ref{thm:mono} establishes the monotonicity part of Theorem~\ref{thm-existence}, and uniqueness is discussed in Theorem~\ref{thm-uni}, and nonexistence of travelling waves for $c < c^*$ is the content of Theorem~\ref{thm-no}.

\section{Notation and auxiliary results}
\label{sec:Notat-auxil-results}

Throughout this paper, $C > 0$ denotes a generic constant, while $C_i$ ($i = 1, 2, \ldots$) represents a specific constant. Let $I$ be an interval, typically $I = \mathbb{R}$. Let $T > 0$ be a real number and $\mathcal{B}$ be a Banach space. We denote by $C([0, T ], \mathcal{B})$ the space of the $\mathcal{B}$-valued continuous functions on $[0, T ]$, while $L^2([0, T ], \mathcal{B})$ is the space of the $\mathcal{B}$-valued $L^2$-functions on $[0, T ]$.  The corresponding spaces of the $\mathcal{B}$-valued functions on $[0,\infty)$ are defined similarly.

For a given travelling wave $\phi$ of (\ref{eq}) satisfying (\ref{eq3}), define
\begin{equation}\label{ch517}
  \begin{split}
    G_j(\xi)&=\partial_jf\left(\phi(\xi),\int_0^{\infty}\int_{\mathbb{R}}h(y,s)
      \phi(\xi-y-cs)dyds\right),\quad     j=1,2,\\
    B(\xi)&=\int_0^{\infty}\int_{\mathbb{R}}h(y,s)G_2(\xi+y+cs)dyds.
  \end{split}
\end{equation}
Obviously, $B(\xi)$ and $G_j(\xi)$, $j=1,2$ are non-increasing and satisfy
\begin{equation}\label{ch517*}
  G_1(\infty)=\partial_1f(K,K),\quad
  B(\infty)=G_2(\infty)=\partial_2f(K,K).
\end{equation}
Moreover, both $G(c,\lambda)$ and $H(\lambda)$ are twice differentiable in $\lambda\in [0, \lambda_0)$. Moreover, $G(c,0)=1$, $H(0)=1$, $H'(\lambda)>0$, $G_{\lambda\lambda}(c,\lambda)>0$, and $H''(\lambda)>0$. Set
\begin{equation}\label{eqde0}
\Delta(c,\lambda)=c\lambda-p\lambda^2-d[H(\lambda)-1]-\partial_1f(0,0)
-\partial_2f(0,0)G(c,\lambda)
\end{equation}
and
\begin{equation}\label{eqde}
  \widetilde{\Delta}(c,\lambda)=c\lambda+p\lambda^2+d[H(-\lambda)-1]
  +\partial_1f(K,K)+\partial_2f(K,K)G(c,-\lambda)
\end{equation}
for all $c\in \mathbb{R}$ and $\lambda\in \mathbb{C}$ with $c\geq 0$ and $\mathrm{Re}\lambda<\lambda^+$, where $\lambda^+=\lambda_0$ if $\partial_2f(0,0)>0$ and $\lambda^+=+\infty$ if $\partial_2f(0,0)=0$.

We require two simple technical statements.
\begin{lemma}\label{lem2}
  There exist $c^*>0$ and $\lambda^*\in (0,\lambda^+)$ such that   $\Delta(c^*,\lambda^*)=0$ and $\Delta_{\lambda}(c^*,\lambda^*)=0$.   Furthermore,
  \begin{description}
  \item[(i)] if $0<c<c^*$, then $\Delta(c,\lambda)<0$ for all     $\lambda\geq 0$;
  \item[(ii)] if $c>c^*$, then the equation $\Delta(c,\cdot)=0$ has     two positive real roots $\lambda_1(c)$ and $\lambda_2(c)$ with     $0<\lambda_1(c)<\lambda^*<\lambda_2(c)<\lambda^+$ such that     $\lambda'_1(c)<0$, $\lambda'_2(c)>0$, $\Delta(c,\lambda)>0$ for     all $\lambda\in (\lambda_1(c),\lambda_2(c))$, and     $\Delta(c,\lambda)<0$ for all $(-\infty,\lambda^+)\setminus     [\lambda_1(c),\lambda_2(c)]$.
\end{description}
\end{lemma}

\noindent \textbf {Proof.} Note that for all $\lambda\in (0,\lambda^+)$,
\begin{align*}
    \Delta_{\lambda}(c,\lambda)&=
    c-2p\lambda-dH'(\lambda)-\partial_2f(0,0)G_{\lambda}(c,\lambda), \\
    \Delta_{\lambda\lambda}(c,\lambda)&=-2p
    -dH''(\lambda)-\partial_2f(0,0)G_{\lambda\lambda}(c,\lambda)<0,
    \\
    \Delta_{c}(c,\lambda)&= \lambda-\partial_2f(0,0)G_c(c,\lambda)>0, \\
    \Delta(c,0)&= -\partial_1f(0,0)
    -\partial_2f(0,0)<0, \\
    \Delta(0,\lambda) &=-2p\lambda^2
    -d[H(\lambda)-1]-\partial_1f(0,0)
    -\partial_2f(0,0)G(0,\lambda)<0
\end{align*}
and
\begin{equation*}
  \lim_{\lambda\to \lambda^+-0}\Delta(c,\lambda)=-\infty.
\end{equation*}
Then the conclusion of this lemma follows.  \hfill $\Box$

\begin{lemma}\label{lem-de2} Under assumptions (F1) and (F2), for each   fixed $c\geq 0$, $\widetilde{\Delta}(c,\cdot)$ has exactly one   positive zero $\upsilon(c)$.
\end{lemma}

\noindent \textbf {Proof.} In view of (F1) and (F2), we have
\begin{equation*}
  \widetilde{\Delta}(c,0)=\partial_1f(K,K)+\partial_2f(K,K)<0
\end{equation*}
and
\begin{equation*}
  \lim_{\lambda\to \lambda^+-0}\widetilde{\Delta}(c,\lambda)=+\infty.
\end{equation*}
Therefore, $\widetilde{\Delta}(c,\cdot)$ has at least one positive zero. Note that
\begin{equation*}
  \widetilde{\Delta}_{\lambda}(c,\lambda)
  =c+2p\lambda-dH'(-\lambda)-\partial_2f(K,K)G_{\lambda}(c,-\lambda)
\end{equation*}
and for $\lambda>0$
\begin{equation*}
  G_{\lambda}(c,-\lambda)=\int^{\infty}_0\int^{\infty}_0h(y,s)\left[
  (y-cs)e^{-\lambda (y+cs)}-(y+cs)e^{\lambda (y-cs)}\right]dyds<0.
\end{equation*}
Then we have $\widetilde{\Delta}_{\lambda}(c,\lambda)>0$ for all $\lambda\in (0,\lambda^+)$. This implies that $\widetilde{\Delta}_{\lambda}(c,\lambda)$ is increasing in $\lambda$ and so it has exactly one positive zero.  \hfill $\Box$

%With the definitions made in Lemmas~\ref{lem2} %and~\ref{lem-de2}, we can introduce
%\begin{equation*}
%  \mathcal{M}(c,\mu)=c\lambda-d[H(\lambda)-1]
%  %-\partial_1f(0,0)-\mu-\partial_2f(0,0)e^{\mu\tau}G(c,\lambda)
%\end{equation*}
%and
%\begin{equation*}
%  \mathcal{N}(\mu)= %\mu+\partial_1f(K,K)+e^{\mu\tau}\partial_2f(K,K),
%\end{equation*}
%where $\lambda$ is any fixed number in $(\lambda_1(c), %\lambda^*)$ if $c>c^*$ and $\lambda=\lambda^*$ if $c=c^*$. We %notice $\mathcal{M}(c,0)=\Delta(c,\lambda)$, which is %positive if $c>c^*$ and is equal to 0 if $c=c^*$. Thus, if %$c>c^*$, there exists $\mu_0>0$ such that %$\mathcal{M}(c,\mu)>0$ for all $\mu\in [0,\mu_0)$. In view of %$\partial_1f(K,K)+\partial_2f(K,K)<0$ (by assumption (F2)), %there exists a small $\mu>0$ so that $\mathcal{N}(\mu)<0$.

We now define the notion of super- and sub-solutions.  For any absolutely continuous function $\varphi\colon\mathbb{R}\to \mathbb{R}$ satisfying that $\varphi'$ and $\varphi''$ exist almost everywhere and are essentially bounded on $\mathbb{R}$, we set
\begin{equation*}
  N_c[\varphi](\xi)\triangleq c\varphi'(\xi)-p\varphi''(\xi)
  -d(J*\phi-\phi)(\xi)-f(\varphi(\xi), (h**\varphi)(\xi)).
\end{equation*}
Given a positive constant $c$, a non-decreasing continuous function $\varphi^+$ is called a \textit{super-solution} of (\ref{eq3}) if $\varphi^+(-\infty)=0$ and $\varphi^+$ is differentiable almost everywhere in $\mathbb{R}$ such that $N_c[\varphi^+](\xi)\geq 0$ for almost every $\xi\in \mathbb{R}$. A continuous function $\varphi^-$ is called a \textit{sub-solution} of (\ref{eq3}) if $\varphi^-(-\infty)=0$, $\varphi^-(\xi)$ is not identically equal to $0$ and $\varphi^-$ is differentiable almost everywhere in $\mathbb{R}$ such that $N_c[\varphi^-](\xi)\leq 0$ for almost every $\xi\in \mathbb{R}$.

Next, we introduce the operator $\mathcal{H}_{\mu} \colon C(\mathbb{R})\to C(\mathbb{R})$ by
\begin{equation*}
  \mathcal{H}_{\mu}(\varphi)=d(J*\varphi-\varphi)
  +f(\varphi, h**\varphi)+\mu \varphi.
\end{equation*}
%for any $\mu>(4d+\max\{|\partial_jf(u,v)|:\, u,v\in [0,K],\, %j=1,2\})/c$.
It is easy to see that $\varphi$ satisfies (\ref{eq3}) if and only if $\varphi$ satisfies
\begin{equation}\label{inte}
  \varphi(\xi)=T_{\mu}(\varphi)(\xi),
\end{equation}
where
\begin{equation*}
  T_{\mu}(\varphi)(\xi)=\frac{1}{c}\int^{\infty}_{0}
  \exp\left\{-\frac{\mu x}{c}\right\}\mathcal{H}_{\mu}(\varphi)(\xi+x)dx
\end{equation*}
if $p=0$, and
\begin{equation*}
  T_{\mu}(\varphi)(\xi)=\frac{1}{p(\varsigma_2-\varsigma_1)}
  \left[\int_{-\infty}^{\xi}
  e^{\varsigma_1(\xi-x)}\mathcal{H}_{\mu}(\varphi)(x)dx+
  \int^{\infty}_{\xi}
  e^{\varsigma_2(\xi-x)}\mathcal{H}_{\mu}(\varphi)(x)dx\right]
\end{equation*}
if $p>0$, and
\begin{equation}\label{varsigma0}
\varsigma_1=\frac{c-\sqrt{c^2+4p\mu}}{2p}<0,\quad
\varsigma_2=\frac{c+\sqrt{c^2+4p\mu}}{2p}>0.
\end{equation}
Choose $\mu>2d+\max\{|\partial_jf(u,v)|:\, (u,v)\in [0,K]\times [0,K],\, j=1,2\}$. Then the operator $T_{\mu}$ is well-defined for functions $\phi$ of a growth rate less than $e^{\mu x}$. Furthermore, since $f$ is monotone in the second argument by (F1), we have for $\varphi \leq \psi$
\begin{equation*}\aligned
  \mathcal{H}_{\mu}(\varphi)- \mathcal{H}_{\mu}(\psi)
  = &\mu [\varphi-\psi]+d[J*(\varphi-\psi)-(\varphi-\psi)]+ \partial_1f(\tilde
  \varphi, h**\varphi) [\varphi - \psi]\leq 0,%\\ &+ \partial_2f(
 % \psi, h**g(\tilde \psi))h**g'(\tilde \psi) [\varphi - \psi]
  \endaligned
\end{equation*}
where $\tilde\varphi(y)$ lies between $\varphi(y)$ and $\psi(y)$.
Then the choice of $\mu$ shows that $\mathcal{H}_{\mu}(\varphi)$ is monotone in $\varphi$,
\begin{equation}\label{mono}
    \mathcal{H}_{\mu}(\varphi)(\xi)\leq \mathcal{H}_{\mu}(\psi)(\xi)
    \quad\mbox{if $0\leq \varphi\leq \psi\leq K$ in $\mathbb{R}$}.
\end{equation}
Thus, we have the following result on the monotonic travelling waves.

\begin{lemma}\label{lemc21} Under assumptions (F1) and (H1), assume that there exists a super-solution   $\varphi^+$ and a sub-solution $\varphi^-$ of (\ref{eq3}) such that   $0\leq \varphi^-\leq \varphi^+\leq K$ on $\mathbb{R}$.
Then (\ref{eq3}) has a solution $\varphi$ satisfying $\varphi'(\xi)\geq   0$ for all $\xi\in \mathbb{R}$.
\end{lemma}

\noindent \textbf {Proof.} Assume that there exist a super-solution $\varphi^+$ and a sub-solution $\varphi^-$ of (\ref{eq3}) such that $0\leq \varphi^-\leq \varphi^+\leq K$ in $\mathbb{R}$.  Define $\varphi_1=T_{\mu}(\varphi^+)$. Then $\varphi_1$ is a well-defined $C^1$ function. From the definition of super-solution, we have
\begin{equation*}
  \varphi^+\geq T_{\mu}(\varphi^+)=\varphi_1.
\end{equation*}
Also, by the definition of sub-solution and the property (\ref{mono}) of $\mathcal{H}_{\mu}$, we get
\begin{equation*}
  \varphi^-\leq T_{\mu}(\varphi^-)\leq T_{\mu}(\varphi^+)=\varphi_1.
\end{equation*}
Hence $\varphi^-(\xi)\leq \varphi_1(\xi)\leq \varphi^+(\xi)$ for all $\xi\in\mathbb{R}$. Moreover, using the fact that $\varphi^+$ is non-decreasing and $\mu>2d+\max\{|\partial_jf(u,v)|:\, u,v\in [0,K],\, j=1,2\}$, we have $\mathcal{H}_{\mu}(\varphi^+)(s)\geq \mathcal{H}_{\mu}(\varphi^+)(\xi)$ for all $s\geq \xi$ and hence
$\varphi'_1(\xi)\geq 0$.
Now define $\varphi_{n+1}=T_{\mu}(\varphi_n)$ for all $n\in \mathbb{N}$. By induction, it is easy to see that $0\leq \varphi^-\leq \varphi_{n+1}\leq \varphi_n\leq \varphi^+\leq K$ and $\varphi'_{n+1}\geq 0$ on $\mathbb{R}$ for all $n\in \mathbb{N}$. Then the limit $\varphi(\xi)\triangleq  \lim_{n\to \infty}\varphi_n(\xi)$ exists for all $\xi\in \mathbb{R}$ and $\varphi(\xi)$ is non-decreasing on $\mathbb{R}$. By Lebesgue's dominated convergence theorem, $\varphi$ satisfies (\ref{inte}) and hence satisfies (\ref{eq3}).

\section{Existence}
\label{sec:Existence}

In this section, we shall establish the existence of travelling waves by constructing a suitable pair of super- and subsolutions. First, we derive two properties of possible solutions of (\ref{eq3}).

\begin{lemma}\label{lem-be} Under assumptions (F1) and (H1), every solution $(c, \varphi)$ of (\ref{eq3}) satisfies $0 <\varphi(\xi)<   K$ for all $\xi\in \mathbb{R}$.
\end{lemma}

\noindent \textbf {Proof.} Let $(c,\varphi)$ be a solution of (\ref{eq3}).  Suppose that there exists $\xi_0\in \mathbb{R}$ such that $\varphi(\xi_0) = 0$. In view of $\varphi(\infty)=K$, without loss of generality, we may assume $\xi_0$ is the right-most point such that $\varphi(\xi_0) = 0$. Since $\varphi(\xi)\geq 0$ for all $\xi\geq \xi_0$, we have $\varphi'(\xi_0)=\varphi''(\xi_0)=0$. It follows from (\ref{eq3}) and (F1) that $\varphi(\xi)\equiv 0$ for all $\xi\in \mathbb{R}$, which contradicts the definition of $\xi_0$. Therefore, $\varphi>0$ on $\mathbb{R}$. Similarly, $\varphi<K$ on $\mathbb{R}$.  This completes the proof. \hfill $\Box$

\begin{lemma}\label{lem-be2}
  Under assumptions (F1) and (H1), every solution $(c, \varphi)$ of   (\ref{eq3}) satisfying $\varphi'\geq 0$ on $\mathbb{R}$ satisfies   $\varphi'>0$ on $\mathbb{R}$.
\end{lemma}

\noindent \textbf {Proof.} Suppose on the contrary that there exists $\xi_0\in \mathbb{R}$ such that $\varphi'(\xi_0) = 0$. By differentiating (\ref{inte}) with respect to $\xi$, we obtain
%\begin{equation*}
%  0=c\varphi'(\xi_0)=\mu e^{\mu\xi_0}\int^{\infty}_{\xi_0}
%  e^{-\mu     %s}\{\mathcal{H}_{\mu}(\varphi)(s)-\mathcal{H}_{\mu}(\varphi)(\xi_0)\}ds
%  \geq 0.
%\end{equation*}
%Thus we have
$\mathcal{H}_{\mu}(\varphi)(s)=\mathcal{H}_{\mu}(\varphi)(\xi_0)$ for all $s\geq \xi_0$. Letting $s\to \infty$, we obtain $\mathcal{H}_{\mu}(\varphi)(\xi_0)=c\mu K$. This, together with $\varphi'(\xi_0)=0$ and (\ref{eq3}), implies that $\varphi(\xi_0)=K$, which contradicts Lemma \ref{lem-be}.  Hence the lemma is proved. \hfill $\Box$

\begin{lemma}\label{lem4}
  Assume that (F1), (F2), (H1) and (H2) hold. Let $c^*$,   $\lambda_1(c)$ and $\lambda_2(c)$ be defined as in Lemma   \ref{lem2}. Let $c > c^*$ be any number. Then for every $\gamma\in (0,\min\{\sigma\lambda_1(c),\lambda_2(c)-\lambda_1(c)\})$ there exists $Q(c,\gamma)>1$ such that  for every $q>Q(c,\gamma)$, the functions $\phi^{\pm}$ defined by
  \begin{equation}\label{eqphi}
    \phi^{+}(\xi)=\min\left\{K,\exp\{\lambda_1(c)\xi\} \right\}
    ,\quad \xi\in \mathbb{R}
  \end{equation}
  and
  \begin{equation}\label{eqphi2}
    \phi^{-}(\xi)= \exp\{\lambda_1(c)\xi\}\max\left\{0,1
      -q\exp\{\gamma\xi\} \right\},
    \quad \xi\in \mathbb{R}
  \end{equation}
  are a super-solution and a sub-solution to (\ref{eq3}),   respectively.
\end{lemma}

\noindent \textbf {Proof.} We begin by proving that $\phi^{\pm}$ are a pair of super- and sub-solutions of (\ref{eq3}).  We only consider the case $\partial_2f(0,0)>0$ because the proof of the case $\partial_2f(0,0)=0$ is similar. It follows from (\ref{eqphi}) that $\phi^{+}(\xi)\leq \exp\{\lambda_1(c)\xi\}$ for all $\xi\in \mathbb{R}$ and hence
$$
J*\phi^{+}(\xi)\leq  \int_{\mathbb{R}}J(y)\exp\{\lambda_1(c)(\xi-y)\}dy= \exp\{\lambda_1(c)\xi\}H(\lambda_1(c))
$$
and
$$
h**\phi^{+}(\xi)\leq  \int^{\infty}_{0}\int_{\mathbb{R}}h(y,s)\exp\{\lambda_1(c)(\xi-y-cs)\}dyds= \exp\{\lambda_1(c)\xi\}G(c,\lambda_1(c)).
$$
Moreover, there exists $\xi^*>0$ satisfying $\exp\{\lambda_1(c)\xi^*\}=K$, $\phi^{+}(\xi)=K$ for $\xi>\xi^*$ and $\phi^+(\xi)=\exp\{\lambda_1(c)\xi\}$ for $\xi\leq \xi^*$.  For $\xi>\xi^*$, we have
\begin{equation*}
  N_c[\phi^{+}](\xi)=-d[J*\phi^{+}(\xi)-K]-
  f(K,(h**\phi^{+})(\xi)
  \geq   -f(K,K)=0.
\end{equation*}
For $\xi\leq \xi^*$, we have
\begin{align*}
    N_c[\phi^+](\xi)
    &\geq
    \phi^{+}(\xi)\{c\lambda_1(c)-p\lambda^2_1(c)-
    d[H(\lambda_1(c))-1]\}
     -f(\phi^{+}(\xi),(h**\phi^{+})(\xi))\\
    &=
    \phi^{+}(\xi)\Delta(c,\lambda_1(c))
    -f(\phi^{+}(\xi),(h**\phi^{+})(\xi))
    +\partial_1f(0,0)\phi^{+}(\xi)\\
    &\quad{}+\phi^{+}(\xi)\partial_2f(0,0)
    G(c,\lambda_1(c))\\
    &\geq
    -f(\phi^{+}(\xi),(h**\phi^{+})(\xi))
    +\partial_1f(0,0)\phi^{+}(\xi)+\partial_2f(0,0)
    (h**\phi^{+})(\xi)\geq 0,
\end{align*}
where we have used the condition (F2) in the last inequality. Therefore, $\phi^{+}$ is a supersolution of (\ref{eq3}).

It follows from (\ref{eqphi2}) that $\exp\{\lambda_1(c)\xi\}\geq \phi^{-}(\xi)\geq \exp\{\lambda_1(c)\xi\}(1
      -q\exp\{\gamma\xi\})$ for all $\xi\in \mathbb{R}$ and hence
$$
\aligned
J*\phi^{-}(\xi) &\geq\int_{\mathbb{R}}J(y)\exp\{\lambda_1(c)(\xi-y)\}(1
      -q\exp\{\gamma(\xi-y)\})dy\\
      & =\exp\{\lambda_1(c)\xi\}H(\lambda_1(c))-q\exp\{[\gamma+\lambda_1(c)]\xi\}H(\gamma+\lambda_1(c)), \\
h**\phi^{-}(\xi)&\geq \int^{\infty}_{0}\int_{\mathbb{R}}h(y,s)\exp\{\lambda_1(c)(\xi-y-cs)\}(1
      -q\exp\{\gamma(\xi-y-cs)\})dyds\\
      &=\exp\{\lambda_1(c)\xi\}G(c,\lambda_1(c))-q\exp\{[\gamma+\lambda_1(c)]\xi\}G(c,\gamma+\lambda_1(c)).
      \endaligned
$$
Let $\xi_0=-\frac{1}{\gamma}\ln q$. Clearly, $\phi^{-}(\xi)=0$ for $\xi>\xi_0$ and $\phi^{-}(\xi)=\exp\{\lambda_1(c)\xi\}(1
      -q\exp\{\gamma\xi\})$ for $\xi\leq \xi_0$. For $\xi>\xi_0$, we have
\begin{equation*}
  N_c[\phi^{-}](\xi)=-dJ*\phi^{-}(\xi)
  -f(0,(h**\phi^{-})(\xi))\leq
  -f(0,0)=0.
\end{equation*}
For $\xi\leq \xi_0$, we have
\begin{align*}
    N_c[\phi^{-}](\xi)
    &\leq     \exp\{\lambda_1(c)\xi\}
    \{c\lambda_1(c)-p\lambda^2_1(c)-d[H(\lambda_1(c))-1]\}\\
    &\quad{}-q\exp\{[\gamma+\lambda_1(c)]\xi\}
    \{c[\gamma+\lambda_1(c)]-p[\gamma+\lambda_1(c)]^2 -d[H(\gamma+\lambda_1(c))-1]\}\\
    &\quad{}-f(\phi^{-}(\xi),(h**\phi^{-})(\xi))\\
    &=  \exp\{\lambda_1(c)\xi\}\Delta(c,\lambda_1(c))
    -q\exp\{[\gamma+\lambda_1(c) ]\xi\}\Delta(c,\gamma+\lambda_1(c))\\
    &\quad{}
    -f(\phi^{-}(\xi),(h**\phi^{-})(\xi))
    +\partial_1f(0,0)\phi^{-}(\xi)+\partial_2f(0,0)\exp\{\lambda_1(c)\xi\}
    G(c,\lambda_1(c))\\
    &\quad{}
    -q\partial_2(0,0)\exp\{[\gamma+\lambda_1(c) ]\xi\}
    G(c,\gamma+\lambda_1(c)) \\
    &\leq
    -q\exp\{[\gamma+\lambda_1(c) ]\xi\}
    \Delta(c,\gamma+\lambda_1(c))
    -f(\phi^{-}(\xi),(h**\phi^{-})(\xi))\\
    &\quad{}+\partial_1f(0,0)
    \phi^{-}(\xi)+\partial_2f(0,0)(h**\phi^{-})(\xi).
\end{align*}
In view of (F2), we have
\begin{align*}
    N_c[\phi^{-}](\xi)
    &\leq
    -q\exp\{[\gamma+\lambda_1(c) ]\xi\}
    \Delta(c,\gamma+\lambda_1(c))+M[\phi^{-}(\xi)
    +(h**\phi^{-})(\xi)]^{1+\sigma}\\
    &\leq
    -q\exp\{[\gamma+\lambda_1(c) ]\xi\}
    \Delta(c,\gamma+\lambda_1(c))+M[1
    +G(c,\gamma+\lambda_1(c))]^{1+\sigma}\exp\{(1+\sigma)\lambda_1(c)\xi\}\\
 &\leq
    \exp\{[\gamma+\lambda_1(c) ]\xi\}
    \Delta(c,\gamma+\lambda_1(c))\left\{\frac{M[1
    +G(c,\gamma+\lambda_1(c))]^{1+\sigma}}{ \Delta(c,\gamma+\lambda_1(c))} -q\right\}
    \leq 0,
\end{align*}
provided that
\begin{equation*}
  q\geq   Q(c,\eta)\triangleq \max\left\{1,\frac{M[1
    +G(c,\gamma+\lambda_1(c))]^{1+\sigma}}{ \Delta(c,\gamma+\lambda_1(c))}\right\}.
\end{equation*}
Therefore, $\phi^{-}$ is a subsolution of (\ref{eq3}). The proof is complete.  \hfill $\Box$

As a consequence of Lemmas \ref{lem-be}, \ref{lem-be2}, and \ref{lem4}, we have the following result on the existence of increasing travelling waves.

\begin{theorem}\label{thm-ex}
Under the conditions (F1), (F2), (H1) and (H2), let $c^*$,
$\lambda_1(c)$ and $\lambda_2(c)$ be defined as in Lemma \ref{lem2}.
Then for each $c> c^*$, (\ref{eq3})
admits a solution $(c, \phi)$ satisfying $0<\phi(\xi)<K$, $\phi'(\xi)>0$ for all $\xi\in \mathbb{R}$, and
\begin{equation}\label{eq4}
\lim\limits_{\xi\to -\infty}\phi(\xi)\exp\{-\lambda\xi\}=1,\quad
\lim\limits_{\xi\to -\infty}\phi'(\xi)\exp\{-\lambda\xi\}=\lambda,
\end{equation} where $\lambda=\lambda_1(c)$ is the smallest positive
zero of $\Delta(c,\cdot)$.
\end{theorem}

\noindent \textbf {Proof.} It follows from Lemmas \ref{lemc21}, \ref{lem-be}, \ref{lem-be2}, and \ref{lem4} that there exists a strictly increasing solution $\phi(\xi)$ to (\ref{eq3}), which will be denoted by $(c,\phi)$ and satisfies
\begin{equation}\label{eq5}
  \exp\{\lambda_1(c)\xi\}\left(1-q\exp\{\gamma\xi\}\right)\leq \phi(\xi)\leq
  \exp\{\lambda_1(c)\xi\},\quad \xi\in
  \mathbb{R}.
\end{equation}
It then follows from (\ref{eq5}) that
\begin{equation*}
  \lim\limits_{\xi\to
    -\infty}\left|\phi(\xi)\exp\{-\lambda_1(c)\xi\}-1\right|\leq
  \lim\limits_{\xi\to -\infty}q\exp\{\gamma\xi\}=0.
\end{equation*}
In view of condition (F2), we have
\begin{align*}
  &\lim\limits_{\xi\to -\infty}\left|f(\phi(\xi),
    (h**\phi)(\xi))-\partial_1f(0,0)\phi(\xi)
    -\partial_2f(0,0)(h**\phi)
    (\xi)\right|\exp\{-\lambda_1(c)\xi\}\\
  &\leq  M\lim\limits_{\xi\to
    -\infty}\left[\phi(\xi)+(h**\phi)(\xi)
  \right]^{1+\sigma}\exp\{-\lambda_1(c)\xi\}\\
  &= M\left[\lim\limits_{\xi\to
      -\infty}\phi(\xi)e^{-\lambda_1(c)\xi/(1+\sigma)}
      +\lim\limits_{\xi\to
      -\infty}(h**\phi)(\xi)\exp\{-\lambda_1(c)\xi/(1+\sigma)\}
  \right]^{1+\sigma}\\
  &= M\left[\lim\limits_{\xi\to
      -\infty}(h**\phi)(\xi)e^{-\lambda_1(c)\xi/(1+\sigma)}
  \right]^{1+\sigma}\\
  &\leq  M\left[g'(0)\lim\limits_{\xi\to
      -\infty}(h**\phi)(\xi)e^{-\lambda_1(c)\xi/(1+\sigma)}
  \right]^{1+\sigma}\\
\end{align*}
and
\begin{align*}
    &\lim\limits_{\xi\to
      -\infty}(h**\phi)(\xi)\exp\{-\lambda_1(c)\xi/(1+\sigma)\}\\
    &=\lim\limits_{\xi\to
      -\infty}\int^{\infty}_{0}\int^{\infty}_{-\infty}
      h(y,t)\phi(\xi-y-ct)
    \exp\{-\lambda_1(c)\xi/(1+\sigma)\}dydt \\
    &= \int^{\infty}_{0}\int^{\infty}_{-\infty}h(y,t)\left[\lim\limits_{\xi\to
        -\infty}\phi(\xi-y-c t)\exp\{-\lambda_1(c)\xi/(1+\sigma)\}\right]dydt= 0.
\end{align*}
Hence, if $p=0$ then for $c>c^*$, we have
\begin{align*}
  &c\lim\limits_{\xi\to -\infty}\phi'(\xi)\exp\{-\lambda_1(c)\xi\}\\
  &= \lim\limits_{\xi\to
    -\infty}\left[d(J*\phi-\phi)(\xi)+f(\phi(\xi),
    (h**\phi)(\xi))\right]\exp\{-\lambda_1(c)\xi\}\\
  &= d[H(\lambda_1(c))-1]+\lim\limits_{\xi\to
    -\infty}f(\phi(\xi),
  (h**\phi)(\xi))\exp\{-\lambda_1(c)\xi\}\\
  &= d[H(\lambda_1(c))-1]+\lim\limits_{\xi\to
    -\infty}\left[\partial_1f(0,0)\phi(\xi)+\partial_2f(0,0)
    (h**\phi)(\xi)\right]\exp\{-\lambda_1(c)\xi\}\\
  &=   d[H(\lambda_1(c))-1]+\partial_1f(0,0)
  +\partial_2f(0,0)G(c,\lambda_1(c))
  = c\lambda_1(c).
\end{align*}
If $p\neq 0$ then for $c>c^*$, using $\phi'(-\infty)=0$ and integrating both sides of (\ref{eq3}) from $-\infty$ to $\xi$, we have
$$
\aligned
&\lim_{\xi\to-\infty}\phi'(\xi)\exp\{-\lambda_1(c)\xi\}\\
=& \frac{c}{p}-\frac{1}{p}\lim_{\xi\to-\infty}\exp\{-\lambda_1(c)\xi\}\int^{\xi}_{-\infty}\left[ d(J*\phi-\phi)(s)+f(\phi(s),
    (h**\phi)(s))\right] ds\\
    =& \frac{c}{p}-\lim_{\xi\to-\infty}\frac{\exp\{-\lambda_1(c)\xi\} \left[ d(J*\phi-\phi)(\xi)+f(\phi(\xi),
    (h**\phi)(\xi))\right]}{p\lambda_1(c)}\\
    =& \frac{c}{p}- \frac{d[H(\lambda_1(c))-1]+\partial_1f(0,0)
  +\partial_2f(0,0)G(c,\lambda_1(c))}{p\lambda_1(c)}
  = \lambda_1(c).
\endaligned
$$
This completes the proof.
\hfill $\Box$

\begin{remark}\label{remark1}
  In Theorem \ref{thm-ex}, by Lebesgue's dominated convergence   theorem, we also have
  \begin{equation*}
    \lim\limits_{\xi\to-\infty}(h**\phi)(\xi)
    \exp\{-\lambda_1(c)\xi\}=G(c,\lambda_1(c)).
  \end{equation*}
\end{remark}

\begin{remark}
  Theorem \ref{thm-ex} implies that the asymptotic behaviours of wave   profiles of the travelling waves obtained by super- and   sub-solutions satisfy (\ref{eq4}) for $c>c^*$. Furthermore, in the   subsequent section, we shall show that the wave profile $\varphi$ of   every travelling wave of (\ref{eq}) satisfying (\ref{eq3}) has   similar asymptotic behaviours.
\end{remark}

Next, we prove that (\ref{eq3}) has a solution $(c, \phi)$ with $0<\phi<K$ and $\phi'>0$ on $\mathbb{R}$ for $c = c^*$.

\begin{theorem}\label{thm-cri}
Under the conditions (F1), (F2), (H1) and (H2),  (\ref{eq3}) has a   solution $(c, \phi)$ with $0<\phi<K$ and $\phi'>0$ on $\mathbb{R}$   for $c = c^*$.
\end{theorem}

\noindent \textbf {Proof.}
We choose a sequence $\{c_j\}\subseteq (c^*,\infty)$ such that $\lim_{j\to\infty}c_j=c^*$. Then for each $j$ there exists a strictly increasing travelling wave $(c_j,\phi_j)$ of (\ref{eq}) such that $\phi_j(-\infty)=0$ and $\phi_j(+\infty)=K$. Since $\phi_j(\cdot+\zeta)$, $\zeta\in \mathbb{R}$, is also a travelling wave, we can assume that $\phi_j(0)=\alpha$ and $\phi_j(x)\leq K$ for a fixed $\alpha\in (0,K)$ and all $x\in \mathbb{R}$ and $j\geq 1$.
Note that $\phi_j$ is a fixed point of operator $T_{\mu}$ in $E$ with $c=c_j$ and
$T_{\mu}(\phi_j)(\xi)$ can be differentiated with respect to $\xi$, where $E$ is the Banach space of bounded and uniformly
continuous functions on $\mathbb{R}$ equipped with the maximum norm. Moreover, we differentiate both sides of (\ref{eq3}) with respect to $\xi$ to
obtain
$$
\aligned
0=&-c\phi''_j(\xi)+p\phi'''_j(\xi)+d[J*\phi'_j-\phi'_j](\xi)  +\partial_1f(\phi_j(\xi),
      (h**\phi_j)(\xi))\phi'_j(\xi)\\
&+\partial_2f(\phi_j(\xi),
      (h**\phi_j)(\xi))h**\phi'_j(\xi).
\endaligned
$$
By the definition of $\mathcal{H}_{\mu}$, it follows that there
exist three positive numbers $N_1$, $N_2$, and $N_3$ (if $p\neq 0$) such that
$$
|\phi'_j(\xi)|\leq N_1,\quad |\phi''_j(\xi)|\leq N_2,\quad |\phi'''_j(\xi)|\leq N_3
$$
for all $n$ and $\xi$. Therefore, $\phi'_j$, $\phi''_j$ and $\phi'''_j$ (if $p\neq 0$) are uniformly bounded and equi-continuous sequences of functions
on $\mathbb{R}$. Then the Arzel\`a-Ascoli theorem implies that there exists
a subsequence of $\{c_j\}$ (for simplicity, denoted again by $\{c_j\}$), such that
$\lim_{j\to\infty}c_j=c^*$, and $\phi'_j$, $\phi''_j$ and $\phi'''_j$ (if $p\neq 0$) converge uniformly on every
bounded and closed subset of $\mathbb{R}$. Thus, $\phi'_j$, $\phi''_j$ and $\phi'''_j$ (if $p\neq 0$) converge
pointwise on $\mathbb{R}$ to $\phi'_*$, $\phi''_*$ and $\phi'''_*$ (if $p\neq 0$), respectively.
By Lebesgue's dominated convergence theorem, letting $j\to \infty$ in the equation $\phi_j=T_{\mu}(\phi_j)$, we then get $\phi_{*}=T_{\mu}(\phi_{*})$.
Thus, $\phi_*$ is a solution of (\ref{eq3}) in the case where $c=c_*$. Clearly, $\phi_*$
is monotonically increasing on $\mathbb{R}$, $\phi_{*}(0)=\alpha$ and $\phi_{*}(x)\leq K$ for all $x\in \mathbb{R}$. One can easily verify that $\phi_{*}(-\infty)=0$ and $\phi_{*}(+\infty)=K$. Thus, (\ref{eq}) has a monotone travelling
wave solution connecting $0$ and $K$ with the wave speed $c=c_*$. This completes the proof.  \hfill $\Box$

\section{Asymptotic behaviour}
\label{sec:Asymptotic-behaviour}

In this section, we always assume that (F1), (F2), (H1) and (H2) hold, and that $c^*$, $\lambda^*$, and $\lambda_1(c)$ are defined as in Lemma \ref{lem2}. We shall follow a method of Carr and Chmaj~\cite{Carr} and Wang, Li, and Ruan~\cite{Wang-Li-Ruan} to establish the exact asymptotic behaviour of the profile $\phi(\xi)$ as $\xi\to -\infty$. For this purpose, we need Ikehara's theorem on the asymptotic behaviour of a positive decreasing function whose Laplace is of a certain given shape.  The proof of Ikehara's theorem can be found, e.g., in~\cite{Carr,Widder}.

\begin{theorem}[Ikehara's theorem]\label{lem-pre1}
  Let $\mathcal{L}[u](\mu)=\int^{\infty}_{0}\exp\{-\mu\xi\}u(\xi)d\xi$ be   the Laplace transform of $u$, with $u$ being a positive   non-decreasing function. Assume that $\mathcal{L}[u]$ has the   representation
  \begin{equation*}
    \mathcal{L}[u](\mu)=\frac{E(\mu)}{(\mu+\alpha)^{k+1}},
  \end{equation*}
  where $k>-1$ and $E$ is analytic in the strip $-\alpha\leq   \mathrm{Re}\mu<0$. Then
  \begin{equation*}
    \lim\limits_{\xi\to\infty}\frac{u(\xi)}{\xi^ke^{-\alpha\xi}}
    =\frac{E(-\alpha)}{\Gamma(\alpha+1)},
  \end{equation*}
  where $\Gamma$ is the Gamma function.
\end{theorem}

\begin{lemma}\label{lem-as1}
  Assume that (F1), (F2), (H1) and (H2) hold. Let $(c, \varphi)$ be a solution of (\ref{eq3}). Then there exists $\gamma>0$ such that
\begin{equation}\label{claim1}
\sup\limits_{\xi\in \mathbb{R}}\varphi(\xi)\exp\{-\gamma\xi\}<\infty\quad \mbox{and}\quad \sup\limits_{\xi\in \mathbb{R}}\Phi(\xi)\exp\{-\gamma\xi\}<\infty,
\end{equation}
where
$\Phi(\xi)\triangleq \int^{\xi}_{-\infty}\varphi(s)ds$.
\end{lemma}

\noindent \textbf {Proof.}
For each $\xi\in \mathbb{R}$, define
$$
\aligned
\psi(\xi)\triangleq & (h**\varphi)(\xi)=
\int_0^{\infty}\int_{\mathbb{R}}h(y,\tau)
\varphi(\xi-y-c\tau)dyd\tau.
\endaligned
 $$
In view of $\partial_2f(0,0)+\partial_1f(0,0)>0$, there exist $\delta_0\in (0,K)$ such that
\begin{equation}\label{assume3}
  \varepsilon_1\triangleq [\partial_2f(0,0)+\partial_1f(0,0)]/4\geq M(u+v)^{\sigma}
\end{equation}
for all $u,v\in [0,\delta_0]$.  In view of (\ref{eqineq}) and $\varphi(-\infty)=0$, there exists $\xi_0<0$ such that for all $\xi<\xi_0$, both $\varphi(\xi)$ and $\psi(\xi)$ lie in the interval $(0,\delta_0)$, where $\delta_0$ is defined as (\ref{assume3}). Thus, for every $\xi<\xi_0$,
\begin{equation}\label{eqineq}
  \begin{split}
    f(\varphi(\xi),\psi(x))
&\geq
    \partial_1f(0,0)\varphi(\xi)+\partial_2f(0,0)\psi(\xi)-M[\varphi(\xi)+\psi(\xi)]^{1+\sigma}\\
&\geq
    \partial_1f(0,0)\varphi(\xi)+\partial_2f(0,0)\psi(\xi)-\varepsilon_1[\varphi(\xi)+\psi(\xi)]\\
    & = (\varepsilon_1-2\varepsilon_2)\varphi(\xi)+(\varepsilon_1+2\varepsilon_2)\psi(\xi)
  \end{split}
\end{equation}
for all $u,v\in [0,\delta_0]$, where $\varepsilon_2 = [\partial_2f(0,0)-\partial_1f(0,0)]/4$. Therefore,
\begin{equation}\label{ineq1}
  c\varphi'(\xi)-p\varphi''(\xi)-d(J*\varphi-\varphi)(\xi)\geq
  (\varepsilon_1-2\varepsilon_2)\varphi(\xi)
  +(\varepsilon_1+2\varepsilon_2)\psi(\xi)
\end{equation}
for all $\xi<\xi_0$.  By using a similar argument as in the proof of Theorem \ref{thm-cri}, we can prove that $\varphi(\xi)$ and $\psi(\xi)$ are both integrable on $(-\infty,0]$.

By Fubini's theorem and Lebesgue's dominated convergence theorem
\begin{align*}
  \int^{\xi}_{-\infty}\psi(s)ds &=
  \int^{\xi}_{-\infty}\left[\int_0^{\infty}\int_{\mathbb{R}}h(y,\tau)
    \varphi(s-y-c\tau)dyd\tau\right]ds\\
  &= \lim\limits_{z\to
    -\infty}\int^{\xi}_{z}\left[\int_0^{\infty}
    \int_{\mathbb{R}}h(y,\tau)
    \varphi(s-y-c\tau)dyd\tau\right]ds\\
  &= \lim\limits_{z\to
    -\infty}\int_0^{\infty}\int_{\mathbb{R}}h(y,\tau)\left[\int^{\xi}_{z}
    \varphi(s-y-c\tau)ds\right]dyd\tau\\
  &=
  \int_0^{\infty}\int_{\mathbb{R}}h(y,\tau)\left[\int^{\xi}_{-\infty}\varphi(s-y-c\tau)ds
  \right]dyd\tau\\
  &= \int_0^{\infty}\int_{\mathbb{R}}h(y)\Phi(\xi-y-c\tau)dyd\tau.
\end{align*}
Integrating (\ref{ineq1}) from $-\infty$ to $\xi$ with $\xi<\xi_0$, we have (using again Fubini's theorem)
\begin{equation}\label{ineq2}
  \begin{split}
    & c\varphi(\xi)-p\varphi'(\xi)-d(J*\Phi-\Phi)(\xi)\\
    &\geq
    (\varepsilon_1-2\varepsilon_2)\Phi(\xi)+(\varepsilon_1+2\varepsilon_2)
    \int_0^{\infty}\int_{\mathbb{R}}h(y,\tau)\Phi(\xi-y-c\tau)dyd\tau\\
    &=
    \varepsilon_1\Phi(\xi)+\varepsilon_1\int_0^{\infty}\int_{\mathbb{R}}h(y,\tau)
    \Phi(\xi-y-c\tau)dyd\tau\\
    &\quad{}+2\varepsilon_2\int_0^{\infty}
    \int_{\mathbb{R}}h(y,\tau)\left[\Phi(\xi-y-c\tau)
      -\Phi(\xi)\right]dyd\tau.
  \end{split}
\end{equation}
Note that
\begin{align*}
  &\int^{\xi}_{z}\int_0^{\infty}\int_{\mathbb{R}}h(y,\tau)
  \left[\Phi(s-y-c\tau)-\Phi(s)
  \right]dyd\tau ds \\
  &= -\int_0^{\infty}\int_{\mathbb{R}}(y+c\tau)h(y,\tau)\int^1_0
  \left[\Phi(\xi-t(y+c\tau))-\Phi(z-t(y+c\tau))\right]dtdyd\tau\\
  &\to
  -\int_0^{\infty}\int_{\mathbb{R}}(y+c\tau)h(y,\tau)
  \int^1_0\Phi(\xi-t(y+c\tau))dtdyd\tau\\
\end{align*}
as $z\to -\infty$. Thus, (\ref{ineq2}) means that $\Phi(\xi)$ and $\int_0^{\infty}\int_{\mathbb{R}}h(y,\tau)\Phi(s-y-c\tau)dyd\tau$ are integrable on $(-\infty,\xi]$.  Moreover,
\begin{multline}\label{ineq3}
  c\Phi(\xi)-p\varphi(\xi)-d\int^{\xi}_{-\infty}(J*\Phi-\Phi)(s)ds
  +2\varepsilon_2\int_0^{\infty}\int_{\mathbb{R}}(y+c\tau)h(y,\tau)
  \int^1_0\Phi(\xi-t(y+c\tau))dtdyd\tau
  \\
  \geq \varepsilon_1\int^{\xi}_{-\infty}\Phi(s)ds
  +\varepsilon_1\int^{\xi}_{-\infty}\int_0^{\infty}
  \int_{\mathbb{R}}h(y,\tau)\Phi(s-y-c\tau)dyd\tau ds.
\end{multline}
Since $\Phi(\xi)$ is increasing, for any $y\in \mathbb{R}$ we have
\begin{equation*}
  (y+c\tau)h(y,\tau)\Phi(\xi)\geq
  (y+c\tau)h(y,\tau)\int^1_0\Phi(\xi-t(y+c\tau))dt.
\end{equation*}
If $\varepsilon_2\geq 0$, then it follows from (\ref{ineq3}) that
\begin{equation}\label{ineq4}
  \aligned
c\Phi(\xi)-p\varphi(\xi)-d\int^{\xi}_{-\infty}(J*\Phi-\Phi)(s)ds
+2\varepsilon_2\Phi(\xi)\int_0^{\infty}
\int_{\mathbb{R}}(y+c\tau)h(y,\tau)dyd\tau
\geq  \varepsilon_1\int^{\xi}_{-\infty}\Phi(s)ds.
\endaligned
\end{equation}
The Mean Value Theorem for Integrals implies that for each $y>0$, there exist $\xi_1(y)\in (\xi,\xi+y)$ and $\xi_2(y)\in (\xi-y,\xi)$ such that $\int^{\xi+y}_{\xi}\Phi(s)ds=y\Phi(\xi_1(y))$ and $\int^{\xi}_{\xi-y}\Phi(s)ds=y\Phi(\xi_2(y))$. It follows from the monotonicity of $\Phi$ that
\begin{equation*}
\aligned
  \int^{\xi}_{-\infty}(J*\Phi-\Phi)(s)ds=&\int^{\xi}_{-\infty}\int_{\mathbb{R}}J(y)[\Phi(x-y)-\Phi(x)]dydx\\
  =&\int_{\mathbb{R}}J(y)\int^{\xi-y}_{\xi}\Phi(x)dxdy\\
  =& \int_0^{\infty}J(y)\int^{\xi-y}_{\xi}\Phi(x)dxdy+\int^0_{-\infty}J(y)\int^{\xi-y}_{\xi}\Phi(x)dxdy\\
  =&  \int_0^{\infty}J(y)\int^{\xi-y}_{\xi}\Phi(x)dxdy+\int_0^{\infty}J(-y)\int^{\xi+y}_{\xi}\Phi(x)dxdy\\
  =& \int_0^{\infty}J(y)\left[\int^{\xi+y}_{\xi}\Phi(x)dx-\int_{\xi-y}^{\xi}\Phi(x)dx\right]dy\\
 =& \int_0^{\infty}yJ(y)[\Phi(\xi_1(y))-\Phi(\xi_2(y))]dy>0.
  \endaligned
\end{equation*}
This, together with (\ref{ineq4}), implies that
\begin{equation}\label{ineq5}
    \left[c+2\varepsilon_2\int_0^{\infty}
    \int_{\mathbb{R}}(y+c\tau)h(y,\tau)dyd\tau\right]\Phi(\xi)
    \geq  \varepsilon_1\int^{\xi}_{-\infty}\Phi(s)ds+p\varphi(\xi).
\end{equation}

If $\varepsilon_2<0$, that is, $\partial_1f(0,0)>\partial_2f(0,0)\geq 0$, then there exists $\xi_0'<\xi_0$ such that
\begin{equation*}
\aligned
  c\varphi'(\xi)-p\varphi''(\xi)-d(J*\varphi-\varphi)(\xi)=& f(\varphi(\xi),\psi(\xi))
  \geq   f(\varphi(\xi),0)\geq
  \frac{1}{2}\partial_1f(0,0)\varphi(\xi)>\varepsilon_1\varphi(\xi)
  \endaligned
\end{equation*}
for all $\xi<\xi'_0$. Thus
\begin{equation}\label{ineq6}
    c\Phi(\xi)
    \geq  \varepsilon_1\int^{\xi}_{-\infty}\Phi(s)ds+p\varphi(\xi).
\end{equation}
Combing (\ref{ineq5}) and (\ref{ineq6}), we have
\begin{multline*}
  \left[c+2|\varepsilon_2|\int_0^{\infty}
  \int_{\mathbb{R}}|y+c\tau|h(y,\tau)dyd\tau\right]\Phi(\xi)
  \geq  \varepsilon_1\int^{0}_{-\infty}\Phi(s+\xi)ds
  \geq   \varepsilon_1\int^0_{-r}\Phi(s+\xi)ds\geq
  \varepsilon_1r\Phi(\xi-r)
\end{multline*}
for all $r>0$ and $\xi<\xi'_0$. Thus there exists $r_0>0$ and some $\theta\in (0,1)$ such that
\begin{equation}\label{gamma}
  \gamma\triangleq \frac{1}{r_0}\ln\frac{1}{\theta}\in (\lambda_1(c),\lambda_0)
\end{equation}
and
\begin{equation*}
  \Phi(\xi-r_0)\leq \theta \Phi(\xi)
\end{equation*}
uniformly in $\xi$. Thus
\begin{equation*}
  \Phi(\xi-r_0)\exp\{-\gamma(\xi-r_0)\}\leq \Phi(\xi)\exp\{-\gamma\xi\}.
\end{equation*}
This, together with $\lim_{\xi\to\infty}\Phi(\xi)e^{-\gamma\xi}=0$, implies that
\begin{equation}\label{claim2}
  \sup\limits_{\xi\in \mathbb{R}}\{\Phi(\xi)\exp\{-\gamma\xi\}\}<\infty.
\end{equation}
Moreover,
$$
\aligned
\exp\{-\gamma\xi\}\int^{\xi}_{-\infty}\Phi(s)ds=&\exp\{-\gamma\xi\}\int^{0}_{-\infty}
\Phi(\xi+s)ds\\
=&\int^{0}_{-\infty}\Phi(\xi+s)\exp\{-\gamma(\xi+s)\}\exp\{\gamma s\}ds
\leq \frac{1}{\gamma}
\sup\limits_{\xi\in \mathbb{R}}\{\Phi(\xi)\exp\{-\gamma\xi\}\}.
\endaligned
$$
Thus, it follows from  (\ref{ineq5}), (\ref{ineq6}), and (\ref{claim2}) that
\begin{equation}\label{claim3}
  \sup\limits_{\xi\in
    \mathbb{R}}\{\varphi(\xi)\exp\{-\gamma\xi\}\}<\infty.
\end{equation}
This completes the proof.  \hfill $\Box$

\begin{lemma}\label{lem-as2}
Assume that (F1), (F2), (H1) and (H2) hold. Let $(c, \varphi)$ be a solution of (\ref{eq3}). Then $\lim_{\xi\to -\infty}\varphi(\xi)\exp\{-\lambda_1(c)\xi\}$ exists for each $c>c^*$.
\end{lemma}

\noindent \textbf {Proof.} Define a bilateral Laplace transform of $\varphi(\xi)$ by
\begin{equation*}
  L[\varphi](\lambda)=\int_{\mathbb{R}}\exp\{-\lambda\xi\}\varphi(\xi)d\xi.
\end{equation*}
By Lemma \ref{lem-as1} and Fubini's theorem, we have
\begin{align*}
  \int_{\mathbb{R}}e^{-\lambda\xi}\psi(\xi)d\xi &=
  \int_{\mathbb{R}}\left[e^{-\lambda\xi}\int_0^{\infty}
  \int_{\mathbb{R}}
    h(y,\tau)\varphi(\xi-y-c\tau)dy\right]d\xi d\tau\\
  &=   \int_0^{\infty}\int_{\mathbb{R}}h(y,\tau)
  \left[\int_{\mathbb{R}}e^{-\lambda\xi}
    \varphi(\xi-y-c\tau)d\xi\right]dy d\tau\\
  &= L(\lambda)\int_0^{\infty}\int_{\mathbb{R}}h(y,\tau)
  e^{-\lambda(y+c\tau)}dyd\tau
  = L(\lambda)G(c,\lambda).
\end{align*}
Take the bilateral Laplace transform of (\ref{eq3}) with respect to $\xi$, we have (with $\Delta(c,\lambda)$ defined in (\ref{eqde0}))
\begin{equation}\label{ineq8}
  \Delta(c,\lambda)L[\varphi](\lambda)=\mathcal{R}(\lambda),
\end{equation}
where $\mathcal{R}(\lambda)$ is the Laplace transform of the function $f(\varphi(\xi),\psi(\xi))-\partial_1f(0,0)\varphi(\xi) -\partial_2f(0,0)\psi(\xi)$. It is not difficult to see that $\mathcal{R}(\lambda)$ is defined for $\lambda$ with $0 < \mathrm{Re}\lambda < \gamma$. In addition, it is easy to see that $\Delta(c,\cdot)$ has no zero $\lambda$ with $\mathrm{Re}\lambda=\lambda_1(c)$ other than $\lambda=\lambda_1(c)$. This implies that $(\lambda-\lambda_1(c))/\Delta(c,\lambda)$ is analytic in the strip $0<\mathrm{Re}\lambda\leq \lambda_1(c)$.  If there exists some $\xi_0>0$ such that $\varphi(\xi)$ is increasing for all $\xi\in (-\infty,-\xi_0)$, then $u(\xi)=\varphi(-\xi)$ is a positive decreasing function on $(\xi_0,\infty)$. Moreover, it follows from (\ref{ineq8}) that
\begin{equation}\label{ineq9}
  \int^{\infty}_{\xi_0}e^{\lambda\xi}u(\xi)d\xi
  =\int_{-\infty}^{-\xi_0}e^{-\lambda\xi}\varphi(\xi)d\xi
  =\frac{\mathcal{E}(\lambda)}{\lambda-\lambda_1(c)}
\end{equation}
with
\begin{equation*}
  \mathcal{E}(\lambda)
  =\frac{(\lambda-\lambda_1(c))\mathcal{R}(\lambda)}{\Delta(c,\lambda)}
  -(\lambda-\lambda_1(c))\int^{\infty}_{-\xi_0}e^{-\lambda\xi}\varphi(\xi)d\xi,
\end{equation*}
which is analytic in the strip $0<\mathrm{Re}\lambda\leq \lambda_1(c)$ because $\int^{\infty}_{-\xi_0}\exp\{-\lambda\xi\}\varphi(\xi)d\xi$ is analytic for all $\mathrm{Re}\lambda>0$.  By means of Theorem \ref{lem-pre1}, $\lim_{\xi\to   \infty}u(\xi)\exp\{\lambda_1(c)\xi\}$, which is equal to $\lim_{\xi\to -\infty}\varphi(\xi)\exp\{-\lambda_1(c)\xi\}$, exists.

If $\varphi(\xi)$ is not monotone on any interval $(-\infty, \xi_0)$ with $|\xi_0|$ sufficiently large, let $\chi(\xi)=\exp\{q\xi\}\varphi(\xi)$, where $q=d/c$ if $p=0$ and
$$
q=\frac{-c+\sqrt{c^2+4pd}}{2p}
$$
if $p> 0$.
Then
\begin{equation*}
  (c+2pq)\chi'(\xi)-p\chi''(\xi)=\exp\{q\xi\}\left[d J*\varphi(\xi)
    +f(\varphi(\xi),\psi(\xi))\right]\geq 0.
\end{equation*}
Suppose that there exists $\xi_1<\xi_2$ such that $\chi(\xi_1)>\chi(\xi_2)$. Note that $\lim_{\xi\to \infty}\chi(\xi)=+\infty$, thus there exists $\xi_3>\xi_1$ such that $\chi'(\xi_3)=0$ and $\chi''(\xi_2)\geq 0$, which contradicts the equation above. Thus, we have $\chi'(\xi)\geq 0$. Then for the bilateral Laplace transform of $\chi(\xi)$, $L[\chi](\lambda)=L(\lambda-q)$. It follows from (\ref{ineq8}) that
\begin{equation*}
  \Delta(c,\lambda-q)L_1(\lambda)=\mathcal{R}(\lambda-q).
\end{equation*}
Using a similar argument as above, we see that
\begin{equation*}
  \lim\limits_{\xi\to
    -\infty}\varphi(\xi)\exp\{-\lambda_1(c)\xi\}=\lim\limits_{\xi\to
    -\infty}\chi(\xi)\exp\{-[q+\lambda_1(c)]\xi\}
\end{equation*}
exists. This completes the proof.  \hfill $\Box$

Using a similar argument as in the proof of the previous lemma, we can verify the following result.
\begin{lemma}\label{lem-as3}
  Assume that (F1), (F2), (H1) and (H2) hold. Let $(c^*, \varphi)$ be   a solution of (\ref{eq3}). Then $\lim_{\xi\to     -\infty}\varphi(\xi)\xi^{-1}\exp\{-\lambda^*\xi\}$ exists.
\end{lemma}

Now we are ready to summarise the asymptotic behaviour of wave profile $\varphi$ as follows.
\begin{theorem}\label{thm-as}
  Under assumptions (F1), (F2), (H1) and (H2), for each solution of   $(c, \varphi)$ of (\ref{eq3}) there exists $\eta=\eta(\varphi)$ such   that
\begin{equation}\label{asy1}
  \lim\limits_{\xi\to-\infty}\frac{\varphi(\xi+\eta)}{\exp\{\lambda_1(c)\xi\}}=1
  \quad\mbox{for}
  \quad c>c^*
\end{equation}
and
\begin{equation}\label{asy2}
  \lim\limits_{\xi\to-\infty}
  \frac{\varphi(\xi+\eta)}{\xi     \exp\{\lambda_1(c)\xi\}}
  =1\quad\mbox{for}
  \quad c=c^*.
\end{equation}
Moreover,
\begin{equation}\label{asy3}
  \lim\limits_{\xi\to-\infty}\frac{\varphi'(\xi)}{\varphi(\xi)}
  =\lambda_1(c)\quad\mbox{for}
  \quad c\geq c^*.
\end{equation}
\end{theorem}

\noindent \textbf {Proof.} Both (\ref{asy1}) and (\ref{asy2}) follow easily from Lemmas \ref{lem-as2} and \ref{lem-as3}. It follows from (\ref{eq3}) that
\begin{align*}
  \lim\limits_{\xi\to-\infty}\frac{c \varphi'(\xi)}{\varphi(\xi)}&=p\lim\limits_{\xi\to-\infty}
  \frac{\varphi''(\xi)}{\varphi(\xi)}+
  d
  \lim\limits_{\xi\to-\infty}\left[\frac{J*\varphi(\xi)}{\varphi(\xi)}
    -1\right]+\lim\limits_{\xi\to-\infty}
  \frac{f(\varphi(\xi),\psi(\xi))}{\varphi(\xi)}\\
  &=p\lambda^2_1(c)+ d\left[H(\lambda_1(c))-1\right]+\partial_1f(0,0)
  +\partial_2f(0,0)\lim\limits_{\xi\to-\infty}\frac{\psi(\xi)}{\varphi(\xi)}\\
  &=p\lambda^2_1(c)+  d\left[H(\lambda_1(c))-1\right]+\partial_1f(0,0) \\
  & \qquad{}
  +\partial_2f(0,0)\lim\limits_{\xi\to-\infty}
  \int_0^{\infty}\int_{\mathbb{R}}h(y,\tau)
  \frac{\varphi(\xi-y-c\tau)}{\varphi(\xi)}dyd\tau\\
  &=p\lambda^2_1(c)+  d\left[H(\lambda_1(c))-1\right]+\partial_1f(0,0)
  +\partial_2f(0,0)G(c,\lambda_1(c))\\
  &= c\lambda_1(c).
\end{align*}
This completes the proof.  \hfill $\Box$

\begin{theorem}\label{thm-as2}
  Under assumptions (F1), (F2), (H1) and (H2), for each solution of   $(c, \varphi)$ of (\ref{eq3}) there exists $\eta=\eta(\varphi)$ such   that
  \begin{equation}\label{asy4}
    \lim\limits_{\xi\to\infty}\frac{K-\varphi(\xi+\eta)}{\exp\{-\upsilon(c)\xi\}}=1,
  \end{equation}
  where $\upsilon(c)$ is the unique positive zero of   $\widetilde{\Delta}(c,\cdot)$, according to Lemma~\ref{lem-de2}.   Moreover,
  \begin{equation}\label{asy5}
    \lim\limits_{\xi\to\infty}\frac{\varphi'(\xi)}{K-\varphi(\xi)}=\upsilon(c).
  \end{equation}
\end{theorem}

\noindent \textbf {Proof.} Define $\Phi(\xi)\triangleq K-\varphi(-\xi)$ and $\Psi(\xi)\triangleq K-(h**\varphi)(-\xi)$.  Obviously, $\Phi(\xi)$ satisfies that $\Phi(-\infty)=0$, $\Phi(\infty)=K$, $0<\Phi(\xi)<K$, and
\begin{equation}\label{ineq10}
  c\Phi'(\xi) = -p\Phi''(\xi)+ f(K-\Phi(\xi),K-\Psi(\xi))-d(J*\Phi-\Phi)(\xi).
\end{equation}
Then for any $\mu>2d+\max\{|\partial_1f(u,v)|:\,u,v\in   [0,K]\}$, we have $\left[\Phi(\xi)e^{-\mu \xi}\right]'<0$ for all $\xi$.  Then, using the bilateral Laplace transform $L[\Phi]$ of $\Phi(\xi)$, we have, using Fubini's theorem again,
\begin{align*}
  \int_{\mathbb{R}}e^{-\lambda\xi}\Psi(\xi)d\xi &=
  \int_{\mathbb{R}}e^{-\lambda\xi}\left[K-\int_0^{\infty}
  \int_{\mathbb{R}}h(y,\tau)
    \varphi(-\xi-y-c\tau)dy d\tau\right]d\xi\\
  &=
  \int_{\mathbb{R}}\left[e^{-\lambda\xi}\int_0^{\infty}
  \int_{\mathbb{R}}h(y,\tau)
    \Phi(\xi+y+c\tau)dyd\tau\right]d\xi\\
  &=
  \int_0^{\infty}\int_{\mathbb{R}}h(y,\tau)\left[\int_{\mathbb{R}}e^{-\lambda\xi}
    \Phi(\xi+y+c\tau)d\xi\right]dyd\tau\\
  &= L[\Phi](\lambda)\int_0^{\infty}\int_{\mathbb{R}}h(y,\tau)
  e^{\lambda(y+c\tau)}dyd\tau
  = L[\Phi](\lambda)G(c,-\lambda).
\end{align*}
Take Laplace transform of (\ref{ineq10}) with respect to $\xi$, we have
\begin{equation}\label{ineq11}
  \widetilde{\Delta}(c,\lambda)\widetilde{L}(\lambda)
  =\widetilde{\mathcal{R}}(\lambda),
\end{equation}
where $\widetilde{\mathcal{R}}(\lambda)$ is the Laplace transform of the function $f(K-\Phi(\xi),K-\Psi(\xi))-\partial_1f(K,K)\Phi(\xi)-\partial_2f(K,K)\Psi(\xi)$. Using a similar arguments as that in the proof of Lemma \ref{lem-as2}, we see that there exists $\eta=\eta(\varphi)$ such that
\begin{equation}\label{ineq12}
  \lim\limits_{\xi\to-\infty}\frac{\Phi(\xi+\eta)}{\exp\{\upsilon(c)\xi\}}=1
  \quad\mbox{and}
  \quad   \lim\limits_{\xi\to-\infty}\frac{\Phi'(\xi)}{\Phi(\xi)}=\lambda_1(c),
\end{equation}
from which (\ref{asy4}) and (\ref{asy5}) follows. The proof is completed.  \hfill $\Box$

\section{Monotonicity and Uniqueness}
\label{sec:Monot-Uniq}

In this section, we investigate the monotonicity and uniqueness (up to a translation) of the travelling wavefront of (\ref{eq3}) by using the sliding method developed in Chen and Guo~\cite{Chen2}.

\begin{theorem}
  \label{thm:mono}
  Under assumptions (F1), (F2), (H1), and (H2), every solution   $(c,\varphi)$ of (\ref{eq3}) satisfies $\varphi'(\xi)>0$ for all   $\xi\in \mathbb{R}$.
\end{theorem}

\noindent \textbf {Proof.} It follows from Lemma \ref{lem-be} and Theorems \ref{thm-as} and \ref{thm-as2} that there exists $M>0$ such that $\varphi'(\xi)>0$ for all $|\xi|\geq M$.  It thus suffices to show that $\varphi'(\xi)>0$ for all $\xi\in [-M,M]$. Suppose on the contrary that $\varphi'(\xi)\leq 0$ for some $\xi_0\in [-M,M]$. By continuity of $\varphi'(\xi)$, there exists $\xi_1\in [-M,M]$ such that $\varphi'(\xi_1)=0$. Then, using the similar arguments as in the proof of Lemma \ref{lem-be2}, we obtain a contradiction. This completes the proof. \hfill $\Box$

In order to prove the uniqueness up to a translation, we shall need the following strong comparison principle.

\begin{lemma}\label{lem5}
Let $(c,\phi_1)$ and $(c,\phi_2)$ be solutions of (\ref{eq3}) with   $\phi_1\geq \phi_2$ on $\mathbb{R}$. Then either $\phi_1\equiv   \phi_2$ or $\phi_1>\phi_2$ on $\mathbb{R}$.
\end{lemma}

\noindent \textbf {Proof.} Suppose that there exists some $\xi_0\in \mathbb{R}$ such that $\phi_1(\xi_0)=\phi_2(\xi_0)$. %Then
%\begin{equation*}
%  0=c\phi_1(\xi_0)-c\phi_2(\xi_0)=\int^{\infty}_0e^{-\mu
%    %x}\left[\mathcal{H}_{\mu}(\phi_1)-\mathcal{H}_{\mu}(\phi_2)\right](\xi_0+x)dx.
%\end{equation*}
In view of $\phi_1\geq \phi_2$ on $\mathbb{R}$, it follows that $\mathcal{H}_{\mu}(\phi_1)(x)=\mathcal{H}_{\mu}(\phi_2)(x)$ for all $x\geq \xi_0$.  It follows from the monotonicity of $\mathcal{H}_{\mu}$ that $\phi_1\equiv \phi_2$ on $\mathbb{R}$.  \hfill $\Box$

\begin{lemma}
  Under assumptions (F1) and (F2), there exists $\varepsilon_0\in   (0,K)$ such that
  \begin{equation}
    \label{ineq50}
    f((1+s)u,(1+s)v)-(1+s)f(u,v)<0
  \end{equation}
  for all $s\in (0,\varepsilon_0)$ and $(u,v)\in \mathbb{R}^2$   satisfying $|K-u|<\varepsilon_0$ and $|K-v|<\varepsilon_0$.
\end{lemma}

\noindent \textbf {Proof.} For each $s\geq 0$, define
\begin{equation*}
  F(s,u,v)=f((1+s)u,(1+s)v)-(1+s)f(u,v).
\end{equation*}
Then $F(0,u,v)=0$ and $F_s(0,u,v)=u\partial_1f(u,v)+v\partial_2f(u,v)-f(u,v)$ for all $(u,v)\in \mathbb{R}^2$.  In view of assumption (F2), we have $F_s(0,K,K)=K\partial_1f(K,K)+K\partial_2f(K,K)<0$. Therefore, there exists $\varepsilon_0>0$ such that $F(s,u,v)<0$ for all $s\in (0,\varepsilon_0)$ and $(u,v)\in \mathbb{R}^2$ satisfying $|K-u|<\varepsilon_0$ and $|K-v|<\varepsilon_0$.  \hfill $\Box$

\begin{lemma}
  \label{lem6}
  Assume that (F1), (F2), (H1) and (H2) hold. Let $(c,\phi_1)$ and   $(c,\phi_2)$ be solutions of (\ref{eq3}). Suppose there exists a   constant $\varepsilon\in (0,\varepsilon_0]$ such that   $(1+\varepsilon)\phi_1(x-\kappa\varepsilon)\geq \phi_2(x)$ on   $\mathbb{R}$, where
  \begin{equation*}
    \kappa=\sup\left\{\frac{\phi_1(x)}{\phi'_1(x)}:\, \phi_1(x)\leq
      K-\varepsilon_0\right\}.
  \end{equation*}
  Then $\phi_1 \geq \phi_2$ on $\mathbb{R}$.
\end{lemma}

\noindent \textbf {Proof.} Define $W(\varepsilon,x)=(1+\varepsilon)\phi_1(x-\kappa\varepsilon)-\phi_2(x)$ and $\varepsilon^*=\inf\{\varepsilon\geq 0$: $W(\varepsilon,x)\geq 0$ for all $x\in \mathbb{R}\}$. By continuity of $W$, $W(\varepsilon^*,x)\geq 0$ for all $x\in \mathbb{R}$. We claim $\varepsilon^*=0$. Suppose on the contrary that $\varepsilon\in (0,\varepsilon_0]$. Then, by the definition of $\kappa$,
\begin{equation*}
  W_{\varepsilon}(\varepsilon,x)=\phi_1(x-\kappa\varepsilon)
  -\kappa(1+\varepsilon)\phi'_1(x-\kappa\varepsilon)<0
\end{equation*}
on $\{x\in \mathbb{R}:\, \phi_1(x-\kappa\varepsilon)\leq K-\varepsilon_0\}$. Noting that $W(\varepsilon^*,\infty)=\varepsilon^*K>0$, we can find $x_0$ with $\phi_1(x_0-\kappa\varepsilon^*)>K-\varepsilon_0$ such that
\begin{equation*}
  0=W(\varepsilon^*,x_0)=W_{x}(\varepsilon^*,x_0)=W_{xx}(\varepsilon^*,x_0).
\end{equation*}
Thus, $(1+\varepsilon^*)\phi_1(\xi_0)=\phi_2(x_0)$, $(1+\varepsilon^*)\phi'_1(\xi_0)=\phi'_2(x_0)$, and $(1+\varepsilon^*)\phi''_1(\xi_0)= \phi''_2(x_0)$,  where $\xi_0=x_0-\kappa\varepsilon^*$. This, together with (\ref{ineq50}), implies
\begin{align*}
  0&=
  -c\phi'_2(x_0)+p\phi''_2(x_0)+d(J*\phi_2-\phi_2)(x_0)+f(\phi_2(x_0),
  (h**\phi_2)(x_0))\\
  &\leq
  -c(1+\varepsilon^*)\phi'_1(\xi_0)+p(1+\varepsilon^*)\phi''_1(\xi_0)+d(1+\varepsilon^*)(J*\phi_1-
  \phi_1)(\xi_0)\\
  &\quad +f((1+\varepsilon^*)\phi_1(\xi_0),\{h**[(1+\varepsilon^*)\phi_1]\}
  (\xi_0))\\
  &=f((1+\varepsilon^*)\phi_1(\xi_0),(1+\varepsilon^*)(h**\phi_1)
  (\xi_0))
   -(1+\varepsilon^*)f(\phi_1(\xi_0),
  (h**\phi_1)(\xi_0))< 0,
\end{align*}
a contradiction. Hence $\varepsilon^*=0$ and so $\phi_1\geq \phi_2$ on $\mathbb{R}$.  \hfill $\Box$

\begin{theorem}\label{thm-uni}
  Assume that (F1), (F2), (H1) and (H2) hold.  For each $c \geq c^*$,   let $(c,\phi_1)$ and $(c,\phi_2)$ be two solutions to (\ref{eq3}).   Then there exists $\gamma\in \mathbb{R}$ such that   $\phi_1(\cdot)=\phi_2(\cdot+\gamma)$, i.e., travelling waves are   unique up to a translation.
\end{theorem}

\noindent \textbf {Proof.} By translating $\phi_2$ if necessary, we can assume that $0<\phi_1(0)=\phi_2(0)<K$. By Theorem \ref{thm-as}, we have
\begin{equation*}
  \lim\limits_{x\to
    -\infty}\frac{\phi_2(x)}{\phi_1(x)}=e^{\lambda_1(c)\theta}
\end{equation*}
for some $\theta\in \mathbb{R}$. Without loss of generality, we assume that $e^{\lambda_1(c)\theta}\leq 1$, for otherwise we can exchange $\phi_1$ and $\phi_2$.  Then
\begin{equation*}
  W(\xi)=\lim\limits_{x\to -\infty}\frac{\phi_2(x)}{\phi_1(x+\xi)}<1
\end{equation*}
for all $\xi>0$. Fix $\xi=1$; then there exists $x_1>0$ such that
\begin{equation}
  \label{for1}
  \mbox{$\phi_1(x+1)>\phi_2(x)$ for all $x\in (-\infty,-x_1)$.}
\end{equation}
Since $\phi_1(\infty)=K$, there exists $x_2\gg 1$ such that $\phi_1(x)\geq K/(1+\varepsilon_0)$ for all $x>x_2$. It follows that
\begin{equation}
  \label{for2}
  \mbox{$(1+\varepsilon_0)\phi_1(x)\geq K\geq \phi_2(x)$ for all     $x>x_2$.}
\end{equation}
Let $\eta=\max\{\phi_2(x):\, x\in [-x_1,x_2]\}\in (0,K)$. In view of $\phi_1(\infty)=K$, there exists $x_3\gg 1$ such that $\phi_1(x)\geq \eta$ for all $x>x_3$. Thus, for $x\in [-x_1,x_2]$, we have $x+x_1+x_3\in [x_3,x_1+x_2+x_3]$ and hence
\begin{equation}
  \label{for3}
  \mbox{$\phi_1(x+x_1+x_3)\geq \eta\geq \phi_2(x)$}.
\end{equation}
Set $z=1+x_1+x_3+\kappa\varepsilon_0$.  It follows from (\ref{for1})--(\ref{for3}) that
\begin{equation}
  \label{for4}
  \mbox{$(1+\varepsilon_0)\phi_1(x+z-\kappa\varepsilon_0)\geq
    \phi_2(x)$ for all $x\in \mathbb{R}$}.
\end{equation}
By monotonicity of $\varphi_1$ and Lemma \ref{lem6}, $\phi_1(x+z)\geq \phi_2(x)$ for all $x\in \mathbb{R}$. Set
\begin{equation*}
  \xi^*=\inf\{z>0:\, \mbox{$\phi_1(x+z)\geq \phi_2(x)$ for all $x\in
    \mathbb{R}$}\}.
\end{equation*}
We claim that $\xi^*=0$. If not, then $\xi^*>0$ and so we have $\phi_1(x+\xi^*)\geq \phi_2(x)$ for all $x\in \mathbb{R}$. It follows from $W(\xi^*/2)<1$ that there exists $x_4>0$ such that
\begin{equation*}
  \mbox{$\phi_1(x+\xi^*/2)\geq \phi_2(x)$ for $x\leq -x_4$}.
\end{equation*}
Consider the function $(1+\varepsilon)\phi_1(x+\xi^*-2\kappa\varepsilon)$. Since $\phi_1(\infty)=K$ and $\phi'_1(\infty)=0$, there exists $x_5\gg 1$ such that
\begin{equation*}
  \frac{d}{d\varepsilon}\left\{(1+\varepsilon)
    \phi_1(x+\xi^*-2\kappa\varepsilon)\right\}=\phi_1(x+\xi^*-2\kappa\varepsilon)-
  2\kappa(1+\varepsilon)\phi'_1(x+\xi^*-2\kappa\varepsilon)>0
\end{equation*}
for all $x\geq x_5$ and $\varepsilon\in [0,1]$. That is, for all $x\geq x_5$ and $\varepsilon\in [0,1]$,
\begin{equation*}
  (1+\varepsilon)\phi_1(x+\xi^*-2\kappa\varepsilon)\geq
  \phi_1(x+\xi^*)\geq \phi_2(x).
\end{equation*}
Now we consider the interval $[-x_4,x_5]$, since $\phi_1(\cdot+\xi^*)\geq \phi_2(\cdot)$, by Lemma \ref{lem5}, $\phi_1(\cdot+\xi^*)>\phi_2(\cdot)$ in $[-x_4,x_5]$. Thus, there exists $\varepsilon\in (0,\min\{\varepsilon_0,\xi^*/(4\kappa)\})$ such that $\phi_1(\cdot+\xi^*-2\kappa\varepsilon)\geq \phi_2(\cdot)$ on $[-x_4,x_5]$. Therefore, combining the estimates on $(-\infty,-x_4]$, $[-x_4,x_5]$, and $[x_5,\infty)$, we conclude that $(1+\varepsilon)\phi_1(\cdot+\xi^*-2\kappa\varepsilon)\geq \phi_2(\cdot)$ on $\mathbb{R}$. It follows from Lemma \ref{lem6} that $\phi_1(\cdot+\xi^*-\kappa\varepsilon)\geq \phi_2(\cdot)$ on $\mathbb{R}$. This contradicts the definition of $\xi^*$. Hence, $\xi^*=0$, i.e., $\phi_1(\cdot)\geq \phi_2(\cdot)$. Since $\phi_1(0)=\phi_2(0)$, we have $\phi_1\equiv \phi_2$ on $\mathbb{R}$. \hfill $\Box$

Finally, we show that no travelling wave solution of speed $c<c^*$ exist. The usual approach is to combine the comparison method
and the finite time-delay approximation to establish the existence of the spreading speed $c^*$ for the solutions with initial
functions having compact supports.  In fact, $c^*$ coincides with the minimal wave speed for monotone travelling waves of
(\ref{eq}). Thus, the nonexistence of travelling waves with the wave speed $c<c^*$ is a straightforward consequence of the
spreading speed. In what follows, however, we shall employ a different method to investigate the nonexistence of travelling waves
with the wave speed $c<c^*$.

\begin{theorem}\label{thm-no}
  Assume that (F1), (F2), (H1) and (H2) hold. Let $c^*$ be defined as in Lemma \ref{lem2}.  Then for every $c\in (0,c^*)$, (\ref{eq}) has no travelling wave front with $(c,\varphi)$ satisfying (\ref{eq3}).
\end{theorem}

\noindent \textbf {Proof.} In view of Theorem \ref{thm:mono}, every solution $(c,\varphi)$ of (\ref{eq3}) satisfies $\varphi'(\xi)>0$ for all $\xi\in \mathbb{R}$.
Take a sequence $\xi_n\to -\infty$ such that $\varphi(\xi_n)\to 0$ and set $v_n(\xi) = \varphi(\xi +\xi_n)/\varphi(\xi_n)$. As $\varphi$ is
bounded and satisfies (\ref{eq3}), the Harnack's inequality implies that the sequence $v_n$ is locally uniformly
bounded. This function $v_n$ satisfies
\begin{equation}
  \label{eq333}
      -cv'_n(\xi)+pv''_n(\xi)+d(J*v_n-v_n)(\xi)+\partial_1f(0,0)v_n(\xi)+\partial_2f(0,0)
      (h**v_n)(\xi)+R_n(\xi)=0
\end{equation}
for $\xi\in \mathbb{R}$,
where
$$
R_n(\xi)=\frac{f(\varphi(\xi+\xi_n),(h**\varphi)(\xi+\xi_n))
-\partial_1f(0,0)\varphi(\xi+\xi_n)-
\partial_2f(0,0)(h**\varphi)(\xi+\xi_n)}{\varphi(\xi_n)}.
$$
The Harnack's inequality implies that the shifted functions $R_n(\xi)$ converge
to zero locally uniformly in $\xi$. Thus one may assume, up to extraction of a subsequence, that the
sequence $v_n$ converges to a function $v$ that satisfies:
\begin{equation}
  \label{eq334}
      -cv'(\xi)+pv''(\xi)+d(J*v-v)(\xi)+\partial_1f(0,0)v(\xi)
      +\partial_2f(0,0)
      (h**v)(\xi)=0,\quad \xi\in \mathbb{R},
\end{equation}
Moreover, $v$ is positive since it is nonnegative and $v(0) = 1$. Equation (\ref{eq334}) admits such a solution
if and only if $c\geq c^*$. Therefore, for every $c\in (0,c^*)$, (\ref{eq}) has no travelling wave front with $(c,\varphi)$ satisfying (\ref{eq3}). This completes the proof.  \hfill $\Box$

%Since $\Delta(c,\lambda)$ has no real zeroes, and $L[\varphi](\lambda)$ is defined for all $\lambda$ satisfying $\mathrm{Re}\lambda> 0$. Thus, (\ref{ineq8}) can be %rewritten as
%$$
%\int_{\mathbb{R}}e^{-\lambda\xi}\left[\Delta(c,\lambda)\varphi(\xi)+\partial_1f(0,0)\varphi(\xi)+\partial_2f(0,0)\psi(\xi)-f(\varphi(\xi),\psi(\xi))\right]d\xi=0
%$$
%This completes the proof.  \hfill $\Box$

\section{Acknowledgements} This work was partially supported by the Natural Science Foundation of People's Republic of China
(Grant No 11671123 \& 11271115), the UK's EPSRC (EP/K027743/1), the Leverhulme Trust (RPG-2013-261) and a Royal Society Wolfson
Research Merit Award.

\end{document}